\documentclass{article}
\usepackage[utf8]{inputenc}
\usepackage{fullpage}
\usepackage{amsmath}
\usepackage{amssymb}
\usepackage{tikz}
\usepackage{svg}
\usepackage{csquotes}
\usepackage{hyperref}
\usepackage{graphicx}
\usepackage{algorithmic}
\usepackage{algorithm}
\usepackage{subfig}
\usepackage[toc,page]{appendix}
\usepackage{amsthm, thmtools}
\usepackage{float}
\usetikzlibrary{positioning}
\usepackage[noabbrev, capitalise]{cleveref}

\theoremstyle{plain}
\newtheorem{prop}{Proposition}

\title{Direct optimization of BPX preconditioners}
\author{Vladimir Fanaskov\thanks{v.fanaskov@skoltech.ru}, Ivan Oseledets\thanks{i.oseledets@skoltech.ru}}

\begin{document}
\maketitle

\begin{abstract}
     We consider an automatic construction of locally optimal preconditioners for positive definite linear systems. To achieve this goal, we introduce a differentiable loss function that does not explicitly include the estimation of minimal eigenvalue. Nevertheless, the resulting optimization problem is equivalent to a direct minimization of the condition number. To demonstrate our approach, we construct a parametric family of modified BPX preconditioners. Namely, we define a set of empirical basis functions for coarse finite element spaces and tune them to achieve better condition number. For considered model equations (that includes Poisson, Helmholtz, Convection-diffusion, Biharmonic, and others), we achieve from two to twenty times smaller condition numbers for symmetric positive definite linear systems.
\end{abstract}

\section{Introduction}
\label{section:Introduction}
In the present contribution, we consider two optimization problems. The first one is the optimization of  a parametric family of preconditioners for a modified Richardson method applied to the matrix $A:A^{T}+A>0$, that is,
\begin{equation}
\label{rho_optimization}
    \omega_{\text{opt}},~\theta_{\text{opt}} = \arg\min_{\omega, \theta} \rho\left(I - \theta B(A, \omega)\right),
\end{equation}
where $B(A, \omega) = B(\omega) A$ (or $B(A, \omega) = B(\omega) A B(\omega)$) is a family of linear systems preconditioned from the left (or in a symmetric fashion), $\rho$ is a spectral radius, and $\omega$ is a set of real numbers. Problem \eqref{rho_optimization} corresponds to a direct optimization of asymptotic convergence speed of an iterative linear solver \cite[Section 2.2.5]{hackbusch1994iterative}.

The second related problem is the optimization of the condition number
\begin{equation}
\label{kappa_optimization}
    \omega_{\text{opt}} = \arg\min_{\omega} \lambda_{\max}\left(B(A, \omega)\right)\big/ \lambda_{\min}\left(B(A, \omega)\right),
\end{equation}
where $\lambda_{\max}$ and $\lambda_{\min}$ are the smallest and the largest eigenvalues, and $A$ is symmetric positive definite.

In both problems we follow the approach adopted in \cite{katrutsa2020black} and further generalized in \cite{greenfeld2019learning}, \cite{luz2020learning}. That is, we introduce a stochastic loss function that approximates an objective function -- spectral radius or a condition number -- and perform a direct gradient-based optimization. The details can be found in \cref{section:Direct_optimization_of_the_spectral_radius} and \cref{section:Direct_optimization_of_the_condition_number}.

For $B(\omega)$ we use a modified BPX \cite{bramble1990parallel} preconditioner. General multilevel preconditioner operates on a chain of linear spaces $V_1 \subset V_2 \subset \dots \subset V_{L}$, where $V_{l},~1\leq l \leq L$ is formed as a linear combination of the set of functions $\phi_{k}^{l}(x),~k=1,\dots, N_{l}$. In the context of a finite element method, $\phi_{k}^{l}(x)$ is a tent function located at vertex $k$ of a grid with the diameter of a cell $\simeq \text{const }2^{-l}$ (grid corresponding to $V_{l+1}$ is constructed from $l$-th grid by, for example, subdivision of coarse triangulation, see i.e. \cite[Section 2]{zhang1992multilevel}). BPX preconditioners were developed for an elliptic problem
\begin{equation}
\label{elliptic_equation}
    -\sum_{i,j=1}^{D}\frac{\partial}{\partial x_{i}}a_{ij}(x) \frac{\partial}{\partial x_{j}}u(x) = f(x),
\end{equation}
with homogeneous Dirichlet boundary conditions and uniformly symmetric positive definite $a_{ij}(x)$. For equation \eqref{elliptic_equation} and a nested set of finite element spaces $\text{span}\left\{\phi_{k}^{l}:~k=1,\dots,N_{l}\right\}$, original BPX and preconditioner reads
\begin{equation}
\label{preconditioners}
    B_{\text{BPX}}(\omega)v = \sum_{l=1}^{L}\sum_{k=1}^{N_{l}}\left(v, \phi_{k}^l\right)\phi_{k}^l,
\end{equation}
where $\left(\psi, \chi\right) = \int \psi(x)\chi(x) dx$ is a $L_2$ scalar product. To improve BPX preconditioner we replace tent function with empirical basis functions $\widetilde{\phi}_{k}^{l},~l=1,\dots, L-1$ and introduce scalars $\widetilde{\alpha}_{l},~l=1,\dots, L-1$ that weight contributions from individual spaces $V_{l}$, that is
\begin{equation}
\label{modified_preconditioners}
    B_{\text{BPX}}(\omega)v = \sum_{l=1}^{L}\widetilde{\alpha}_{l}\left(\sum_{k}\left(v, \widetilde{\phi}_{k}^l\right)\widetilde{\phi}_{k}^l\right).
\end{equation}
The details of the parametrisation and more convenient form of preconditioners \eqref{modified_preconditioners} are given in \cref{section:BPX_precs}.

Together $\widetilde{\phi}_{k}^{l}$ and $\widetilde{\alpha}_{l}$ form a set of parameters $\omega$ in problems \eqref{rho_optimization}, \eqref{kappa_optimization}. The results of the optimization can be found in \cref{section:Numerical_examples}. In short, our framework allows for up to two times smaller spectral radius of modified Richardson scheme and up to twenty times smaller condition number for selected problems.

\section{Direct optimization of the spectral radius}
\label{section:Direct_optimization_of_the_spectral_radius}

Problem \eqref{rho_optimization} can be viewed in the context of a general search for better linear iterative methods. As explained in \cite[Section 2.2.2]{hackbusch1994iterative}, an arbitrary consistent iterative method can be written in a form
\begin{equation}
\label{linear_iteration}
    x^{n+1} = M(\omega, A)x^{n} + N(\omega, A)b,~I-M(\omega, A) = N(\omega, A)A.
\end{equation}
The efficiency of the method can be characterised by spectral radius $\rho\left(M(\omega, A)\right)$, because it quantifies an asymptotic convergence rate in a following sense. Let $e^{n}$ be an error vector on step $n$, $\left\|\cdot\right\|$ is arbitrary norm and $\rho_{m+k, m} = \left(\left\|e^{m+k}\right\|\big/\left\|e^{m}\right\|\right)^{1\big/k}$ is a geometric mean of a one-step error reduction factor $\rho_{m+1, m}$. It is known that $\lim_{k\rightarrow \infty} \max_{x_0}\left\{\rho_{m+k, m}(x_0)\right\} = \rho(M(\omega, A))$ (see \cite[Remark 2.22]{hackbusch1994iterative}). That is, $\rho\left(M(\omega, A)\right)$ characterises a geometric mean of an error reduction per iteration in the worst case. Because of that it is a custom to use $\rho\left(M(\omega, A)\right)$ as an objective function. For example, classical schemes like SOR and instationary Richardson iteration were optimized analytically \cite{hadjidimos2000successive}, \cite[chapters 4, 8]{hackbusch1994iterative}  and numerically \cite{manteuffel1978adaptive}, \cite{reid1966method}, to achieve better $\rho\left(M(\omega, A)\right)$. More modern attempts include optimization of multigrid with local Fourier analysis \cite{brown2020tuning} and directly \cite{schmitt2019optimizing}, \cite{luz2020learning}, \cite{greenfeld2019learning}, \cite{katrutsa2020black}.

To apply gradient-based optimization to \eqref{rho_optimization} we need a differentiable approximation to the spectral radius. We consider three options.

The first one is an approximation of $\rho(A)$ by Gelfand formula \cite{kozyakin2009accuracy} $\rho(A)= \lim_{k\rightarrow\infty}\left\|A^{k}\right\|^{1/k}$ combined with a stochastic trace approximation \cite{avron2011randomized}:
\begin{equation}
    \label{stochastic_trace}
    \begin{split}
    \rho(A) \simeq \rho_1(A, k, N_{\text{batch}}) \equiv \left(\frac{1}{N_{\text{batch}}}\sum_{j=1}^{N_{\text{batch}}} \left\|A^{k} z_j\right\|_{2}^{2}\right)^{1\big/2k},\\
    \forall j:\mathbb{P}\left((z_{j})_{i}=\pm 1\right) = 1\big/2,~\forall i, j:~z_{i},z_{j} \text{ are independent}.
    \end{split}
\end{equation}
More details about this approach can be found in \cite{katrutsa2020black}.

The second option is based on $\rho(A) = \lim_{k\rightarrow \infty} \left(\left\|e^{m+k}\right\|\big/\left\|e^{m}\right\|\right)^{1\big/k},~e^{m+l} = A^{l} e^{m}$, see \cite[Remark 2.22 (b)]{hackbusch1994iterative} for details. This gives us another approximation
\begin{equation}
    \label{error_propagation}
    \rho(A) \simeq \rho_2\left(A, k\right) \equiv \left( \left\|A^{k} z\right\|_{2}\big/\left\|z\right\|_2\right)^{1\big/k},~\left(z_{i}\right)_{j} \sim \mathcal{N}(0, 1).
\end{equation}
Approximation \eqref{error_propagation} does not contain averaging, but we can introduce $N_{\text{batch}}$ the same way as in \eqref{stochastic_trace}. That gives us the following the last approximation
\begin{equation}
    \label{error_propagation_2}
    \begin{split}
    &\rho(A) \simeq \rho_3\left(A, k, N_{\text{batch}}\right) \equiv \frac{1}{N_{\text{batch}}}\sum_{j=1}^{N_{\text{batch}}}\left( \left\|A^{k} z_{j}\right\|_{2}\big/\left\|z_{j}\right\|_2\right)^{1\big/k},\\
    &\forall j:\left(z_{j}\right)_{i} \sim \mathcal{N}(0, 1),~\forall i, j:~z_{i},z_{j} \text{ are independent}.
    \end{split}
\end{equation}
The resulting loss will measure how well matrix $A$ damps nonzero initial vectors on average. We observed that introduction of $N_{\text{batch}}>1$ in \eqref{error_propagation_2} leads to better convergence.

With approximations $\rho_i\left(A, k, N_{\text{batch}}\right),~i=1,2,3$ we can use forward mode automatic differentiation \cite{RevelsLubinPapamarkou2016} and standard optimizers \cite[Section 8.3]{goodfellow2016deep}  to solve problem \eqref{rho_optimization}. The resulting algorithm coincides with \cref{algorithm:rho_optimization} with $N_{\text{inner}}=1$.

\section{Direct optimization of the condition number}
\label{section:Direct_optimization_of_the_condition_number}

\begin{figure}[htbp]
  \centering
  \includegraphics[scale=0.55]{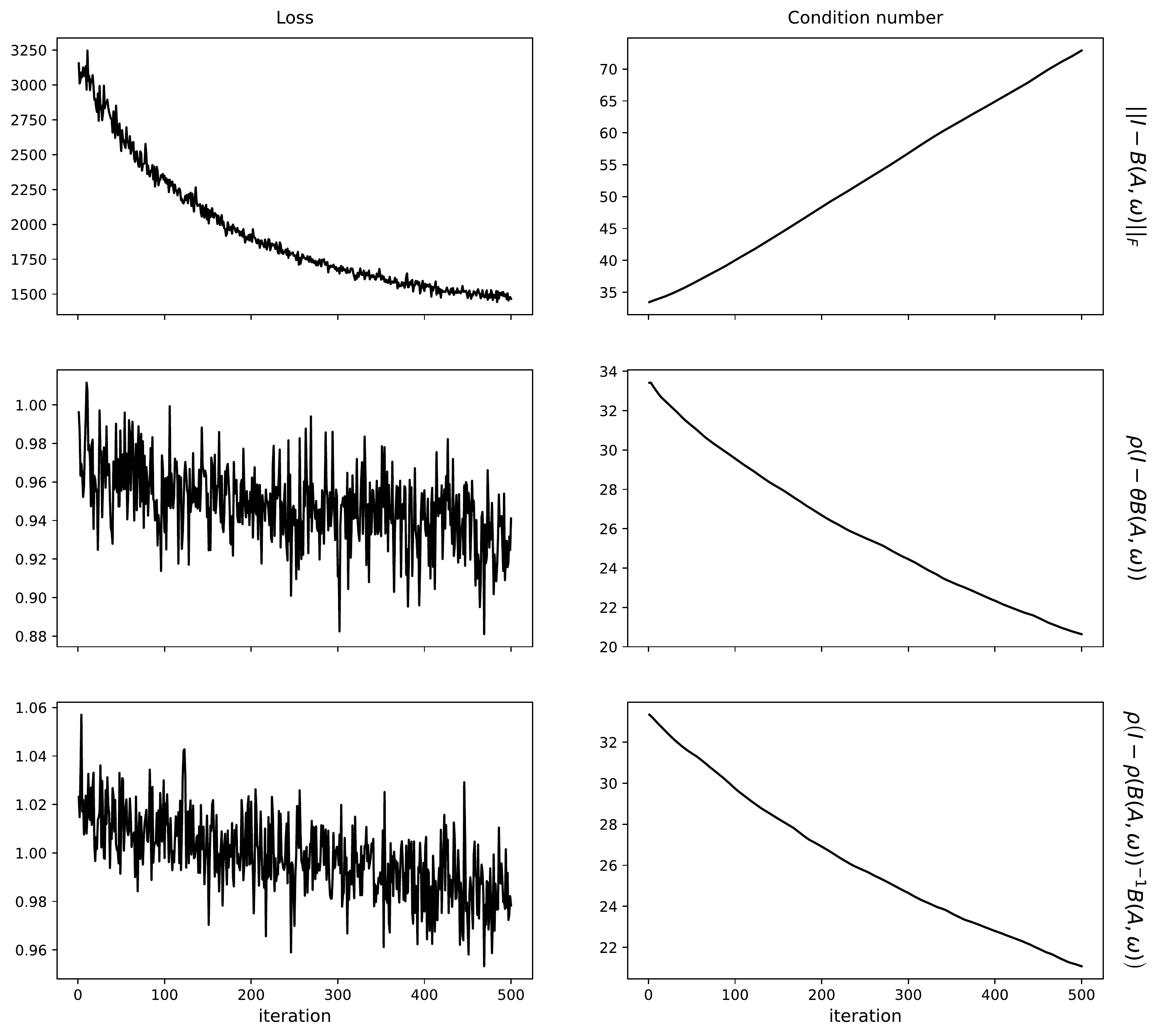}
  \caption{Comparison of three loss functions. The first column shows how the value of the loss function changes in the course of iterations, graphs in the second column demonstrate an evolution of condition number. The first row corresponds to the Frobenius norm $\left\|I - B(\omega)A\right\|$ used as a loss function, the second row shows minimization of $L_1$ by \cref{algorithm:rho_optimization} ($N_\text{inner}=1$), the last row shows minimization of $L_2$ by \cref{algorithm:kappa_optimization}. For the last two cases, we used \eqref{stochastic_trace} to approximate spectral radius. It is clear that the decrease of both losses $L_1$ and $L_2$ lead to a smaller spectral condition number, whereas smaller Frobenius norm does not lead to a better spectral condition number. In all cases we use modified BPX preconditioner \eqref{modified_BPX} as $B(\omega)$ and FEM discretization (see \cref{section:BPX_precs}) of Poisson equation \eqref{Poisson_equation} in $D=1$.}
  \label{fig:compare_losses}
\end{figure}

Unlike problem \eqref{rho_optimization} the optimization of the condition number is not straightforward. The main problem is the presence of $\lambda_{\min}$ which is not readily available. The standard way to resolve this issue is to substitute spectral radius with more amenable loss. For example, objective functions $\left\|R - A\right\|$ and $\left\|I - R^{-1}A\right\|$ (here $R$ is an easy invertible approximation to $A$) were used to construct optimal circulant \cite{chan1988optimal}, \cite{tyrtyshnikov1992optimal}, \cite{strang1986proposal} and sparse approximate inverse \cite{grote1997parallel}, \cite{chow1998approximate} preconditioners. It is known that for nonsymmetric matrices optimization of $\left\|I - R^{-1}A\right\|$ can fail to deliver good preconditioner \cite{chow1994approximate}. The same is true for symmetric positive definite matrices as illustrated on \cref{fig:compare_losses}.

For symmetric positive definite matrices, one can construct a loss function that leads to a direct minimization of the spectral condition number. It is well known that for arbitrary positive definite matrix $C$, optimal spectral radius of $I-\theta C$ is $\left(\lambda_{\max}(C)-\lambda_{\min}(C)\right)\big/\left(\lambda_{\max}(C)+\lambda_{\min}(C)\right)$. Using this fact, we can consider the following loss function
\begin{equation}
\label{kappa_optimization_2}
    L_1(\omega) = \rho\left(I - \theta_{\text{opt}}(\omega) B(A, \omega)\right),~\theta_{\text{opt}}(\omega) = \arg\min_{\theta} \rho\left(I - \theta B(A, \omega)\right).
\end{equation}

\begin{algorithm}[t]
    \caption{Minimization of $L_1$ \eqref{kappa_optimization_2}.}
    \label{algorithm:rho_optimization}
    \begin{algorithmic}
        \STATE Input: matrix $A>0$, parametric family of preconditioners $B(A,\omega): B(A,\omega)>0 $, stochastic gradient-based optimizer $\omega \leftarrow O\left(\omega, \partial_{\omega}\left(\text{loss function}\right)\right)$ (f.e., ADAM, \cite{kingma2014adam}), batch size $N_{\text{batch}}$, number of matrix-vector products $k$, number of epochs $N_{\text{epochs}}$, number of iterations for inner loop $N_{\text{inner}}$, estimator of the spectral radius $m\in\left\{1, 2, 3\right\}$.
        \\\hrulefill
        \FOR{$i = 1:N_{\text{epochs}}$}
            \FOR{$j = 1:N_{\text{inner}}$}
                \STATE $\rho_{m}, \partial_{\theta} \rho_{m} \leftarrow \text{AD }\rho_{m}\left(I-\theta B(A, \omega), k, N_{\text{batch}}\right)$ // AD -- automatic differentiation
                \STATE $\theta \leftarrow O\left(\theta, \partial_{\theta} \rho_{m} \right)$
            \ENDFOR
            \STATE $L_1, \partial_{\omega} L_1 \leftarrow \text{AD } \rho_{m}\left(I-\theta B(A, \omega), k, N_{\text{batch}}\right)$
            \STATE $\omega \leftarrow O\left(\omega, \partial_{\omega} L_1\right)$
        \ENDFOR
    \end{algorithmic}
\end{algorithm}

Evidently, the minimization of \eqref{kappa_optimization_2} is equivalent to the minimization of $\left(\kappa(B(A, \omega)) - 1\right)\big/\left(\kappa(B(A, \omega)) + 1\right)$, where $\kappa$ is the spectral condition number. That means we constructed an optimization problem equivalent to \eqref{kappa_optimization} but without $\lambda_{\min}$. A procedure for minimization of loss \eqref{kappa_optimization_2} is summarised in \cref{algorithm:rho_optimization}. The inner loop finds $\theta_{\text{opt}}$ for each $\omega$ and the outer loop optimizes $\omega$. If an inner loop is reduced to a single iteration as it is done in many other situations (for example, generalized policy iteration \cite[Section 4.6]{sutton1998introduction}, and full approximation scheme \cite[Section 5.3.1]{trottenberg2000multigrid} follow the same pattern), we obtain an algorithm that minimizes spectral radius for modified Richardson scheme.

Another equivalent loss function is
\begin{equation}
\label{kappa_optimization_3}
    L_2(\omega) = \rho\left(I - \frac{1}{\rho(B(A, \omega))}B(A, \omega)\right).
\end{equation}
Indeed, $\rho\left(I - B(A, \omega)/\rho(B(A, \omega))\right) = 1 - \lambda_{\min}(B(A, \omega))\big/\lambda_{\max}(B(A, \omega))$, which means that a minimization of \eqref{kappa_optimization_3} is equivalent to minimization of $1-1\big/\kappa(B(A, \omega))$. Gradient-based optimization can be applied to \eqref{kappa_optimization_3} directly, but we can exploit a special structure of the problem to shorten the computation graph. Using a chain rule we get
\begin{equation}
\label{chain_rule}
\begin{split}
    &\frac{\partial}{\partial \omega_{i}} \rho\left(I - \frac{1}{\rho(B(A, \omega))}B(A, \omega)\right) = \left.\left(\frac{\partial}{\partial \omega_{i}} \rho\left(I - \theta B(A, \omega)\right)\right)\right|_{\theta=\rho(B(A, \omega))^{-1}} \\
    &- \left.\left(\theta^2\frac{\partial}{\partial \theta} \rho\left(I - \theta B(A, \omega)\right)\right)\right|_{\theta=\rho(B(A, \omega))^{-1}}\frac{\partial }{\partial \omega_{i}}\rho(B(A, \omega)).
    \end{split}
\end{equation}
This leads to \cref{algorithm:kappa_optimization}. The performance of these two loss function is illustrated on \cref{fig:compare_losses}. In our experiments, we find little difference between \cref{algorithm:rho_optimization} and \cref{algorithm:kappa_optimization}. Because of that, we mainly use \cref{algorithm:rho_optimization}, which requires a single computation of a gradient with respect to $\omega$. However, unlike $L_1$ loss function $L_2$ is defined in terms of $\rho$ in closed form, i.e., without an additional optimization problem, so it can be more advantageous in situations when a family of preconditioners is learned for a set of related linear equations, as it is done in \cite{greenfeld2019learning} for the multigrid solver.

\begin{algorithm}[t]
    \caption{Minimization of $L_2$ \eqref{kappa_optimization_3}.}
    \label{algorithm:kappa_optimization}
    \begin{algorithmic}
        \STATE Input: symmetric positive definite matrix $A>0$, parametric family of preconditioners $B(A,\omega): B(A,\omega)>0 $, stochastic gradient-based optimizer $\omega \leftarrow O\left(\omega, \partial_{\omega}\left(\text{loss function}\right)\right)$ (f.e., ADAM, \cite{kingma2014adam}), batch size $N_{\text{batch}}$, number of matrix-vector products $k$, number of epochs $N_{\text{epochs}}$, estimator of the spectral radius $m\in\left\{1, 2, 3\right\}$.
        \\\hrulefill
        \FOR{$i = 1:N_{\text{epochs}}$}
            \STATE $\theta \leftarrow 1\big/\rho_{1}\left(B(A, \omega), k, N_{\text{batch}}\right)$
            \STATE $\rho_{m}, \partial_{\omega} \rho_{m} \leftarrow \text{AD } \rho_{m}\left(B(A, \omega), k, N_{\text{batch}}\right)$ // AD -- automatic differentiation
            \STATE $L_2, \partial_{\omega} L_2, \partial_{\theta} L_2 \leftarrow \text{AD } \rho_{m}\left(I-\theta B(A, \omega), k, N_{\text{batch}}\right)$
            \STATE $\omega \leftarrow O\left(\omega, \partial_{\omega} \rho_{m}-\theta^2 \partial_{\theta} L_2 \partial_{\omega} L_2\right)$
        \ENDFOR
    \end{algorithmic}
\end{algorithm}

We summarize the results of this section in the following statement.

\begin{prop}
\label{prop:optimization}
    Let $A>0$ and $B(\omega)>0$ for all $\omega$. For left $B(A, \omega) = B(\omega)A$, symmetric $B(A, \omega) = B(\omega) A B(\omega)$ and  right $B(A, \omega) = A B(\omega)$ preconditioners the following three optimization problems are equivalent:
    \begin{itemize}
        \item $\min_{\omega} \rho\left(I - \theta_{\text{opt}}(\omega) B(A, \omega)\right)$, where $\theta_{\text{opt}}(\omega) = \arg\min_{\theta} \rho\left(I - \theta B(A, \omega)\right)$ -- loss function \eqref{kappa_optimization_2}
        \item $\min_{\omega} \rho\left(I - B(A, \omega)\big/ \rho(B(A, \omega))\right)$ -- loss function \eqref{kappa_optimization_3}
        \item $\min_{\omega} \left(\lambda_{\max}(B(A, \omega)) \big/ \lambda_{\min}(B(A, \omega))\right)$
    \end{itemize}
\end{prop}

\section{Modified BPX preconditioner}
\label{section:BPX_precs}

\begin{figure}
    \centering
    \subfloat[][Bilinear FEM]{
        \begin{tabular}{c||c|c|c||c|c|c}
     & \multicolumn{3}{|c||}{BPX} & \multicolumn{3}{|c}{optimized BPX}\\
    \hline
    $L$ & $\rho$ & $\kappa$& $N$ &  $\rho$ & $\kappa$ & $N$\\
    \hline
    $3$ & $0.621$ & $4.277$ & $5$ & $0.314$ & $1.915$ & $2$\\
    \hline
    $4$ & $0.701$ & $5.678$ & $7$ & $0.386$ & $2.259$ & $3$\\
    \hline
    $5$ & $0.746$ & $6.867$ & $8$ & $0.432$ & $2.523$ & $3$\\
    \hline
    $6$ & $0.774$ & $7.866$ & $10$ & $0.46$ & $2.706$ & $3$\\
\end{tabular}
        \label{fig:Poisson_2D:FEM}
    }
    \subfloat[][Mehrstellen]{
        \begin{tabular}{c||c|c|c||c|c|c}
     & \multicolumn{3}{|c||}{BPX} & \multicolumn{3}{|c}{optimized BPX}\\
    \hline
    $L$ & $\rho$ & $\kappa$& $N$ &  $\rho$ & $\kappa$ & $N$\\
    \hline
    $3$ & $0.62$ & $4.269$ & $5$ & $0.427$ & $2.488$ & $3$\\
    \hline
    $4$ & $0.7$ & $5.678$ & $7$ & $0.448$ & $2.621$ & $3$\\
    \hline
    $5$ & $0.746$ & $6.867$ & $8$ & $0.454$ & $2.666$ & $3$\\
    \hline
    $6$ & $0.774$ & $7.867$ & $10$ & $0.472$ & $2.791$ & $4$\\
\end{tabular}
        \label{fig:Poisson_2D:Mehrstellen}
    }
    \caption{Results of optimization for $2\text{D}$ Poisson equation \eqref{Poisson_equation}. Here $\rho = \lambda_{\max}\left(I-\theta_{\text{opt}}BAB\right)$ -- a spectral radius of optimal Richardson iteration for a given preconditioner, $\kappa = \lambda_{\max}(BAB)\big/\lambda_{\min}(BAB)$ -- spectral condition number, and $N$ -- the number of iteration needed to drop an error by $0.1$ in the arbitrary norm, i.e., $\left\|e^{n+N}\right\|\big/\left\|e^{n}\right\| \leq 0.1$.}
    \label{fig:Poisson_2D}
\end{figure}

We already specified algorithms that can be used to optimize condition number (optimization problem \eqref{kappa_optimization}). In this section, we describe a parametric family of positive definite preconditioners that we use in optimization.

To obtain a convenient form of BPX preconditioner, we introduce a hierarchy of meshes
\begin{equation}
\label{hierarchy_of_meshes}
    M_{l} = \left\{x_{j}^{l} = j\big/2^{l}:j=0,1,\dots~,2^l-1,2^l\right\},~l=1,\dots,L
\end{equation}
such that each next mesh contains a previous one, that is, $M_{l} \subset M_{l+1}$. For each mesh, we define a set of basis functions $\phi_{i}^{l}(x) = \phi^{l}(x-x_{i}),i=0,\dots,2^{l}$, which are rescaled and translated copies of a tent function $\phi^{l}(x) = \left(1 + x\big/2^l\right)\text{Ind}\left[-1\big/2^{l}\leq x\leq0\right] + \left(1 - x\big/2^l\right)\text{Ind}\left[0<x\leq1\big/2^{l}\right]$, where $\text{Ind}\left[x\right]$ is $1$ if $x$ holds and $0$ otherwise. Basis functions $\left\{\phi^{L}_{i}(x):i=0,\dots,2^{L}\right\}$ are used to perform standard finite element discretization \cite{ciarlet2002finite} of elliptic problem \eqref{elliptic_equation} for $x\in\left[0, 1\right]$. For higher dimensions, we use $M_{l}$ and $\phi^{l}_{i}$ that are direct products of unidimensional meshes and basis functions.

In article \cite{bachmayr2020stability}, authors show that for equation \eqref{elliptic_equation} in $D=1$ with uniform Dirichlet boundary condition at $x=0$ and uniform Neumann boundary condition at $x=1$ discretized as we just described, BPX preconditioner has the following form
\begin{equation}
\label{BPX_matrix_form}
    \begin{split}
        &\mathcal{B} = \sum_{k=1}^{L}\alpha_{k}B_{k}^{L}B_{L}^{k},~B_{l}^{L} = I_{l}\otimes\eta_{L-l} + S_{l}\otimes\left(\xi_{L-l}-\eta_{L-l}\right), B_{L}^{l} = \left(B_{l}^{L}\right)^T,~\alpha_{k}=1\\
        &\left(\eta_{k}\right)_{i}=i/2^{k},\left(\xi_{k}\right)_i = 1,~\left(S_{l}\right)_{ij} = \delta_{ij+1},\left(I_{l}\right)_{ij} =\delta_{ij},~i,j=1,~\dots~,~2^{l}.
    \end{split}
\end{equation}
If $D=2$ matrices $B_{L}^{k}$ are replaced with $B_{L}^{k}\otimes B_{L}^{k}$ and $\alpha_{k}$ are with ratio of grid spacings $h_{L}\big/h_{k}$. The proof of the optimality of symmetric preconditioner \eqref{BPX_matrix_form} can be found in \cite[Appendix A]{bachmayr2020stability}.

It is easy to see that components of $\eta_{L-l}$ and $\xi_{L-l} - \eta_{L-l}$ contains scalar products $\left(\phi^{L}, \phi^{l}\right)$. Using this observation, one can extend \eqref{BPX_matrix_form} on other boundary conditions:
\begin{prop}
\label{prop:BPX_BCs}
    For equation \eqref{elliptic_equation} in $D=1$ discretized with linear finite elements, symmetric BPX preconditioner has a form $\mathcal{B} = \sum_{k=1}^{L}\alpha_{k}B_{k}^{L}B_{L}^{k}$, where matrices $B_{L}^{k}$ depend on boundary conditions as follows:
    \begin{itemize}
        \item Dirichlet-Neumann: $B_{l}^{L} = I_{l}\otimes\eta_{L-l} + S_{l}\otimes\left(\xi_{L-l}-\eta_{L-l}\right);$
        \item Neumann-Dirichlet: $B_{l}^{L} = I_{l}\otimes\eta^{r}_{L-l} + \left(S_{l}\right)^{T}\otimes\left(\xi_{L-l}-\eta^{r}_{L-l}\right),~\left(\eta^{r}_{k}\right)_{i} = \left(\eta_{k}\right)_{2^{k}-i+1};$
        \item Neumann-Neumann: $B_{l}^{L} = \begin{pmatrix} 1&0_{1\times2^{l}}\\ e_l\otimes \left(\xi_{L-l}-\eta_{L-l}\right)&I_{l}\otimes\eta_{L-l} + S_{l}\otimes\left(\xi_{L-l}-\eta_{L-l}\right) \end{pmatrix};$
        \item Dirichlet-Dirichlet: $B_{l}^{L} = \left[I_{l}\otimes\eta_{L-l} + S_{l}\otimes\left(\xi_{L-l}-\eta_{L-l}\right)\right]_{\text{last row and column are removed}}.$
    \end{itemize}
     All boundary conditions are uniform and vectors $\xi_{L-l}, \eta_{L-l}$ are defined as in \eqref{BPX_matrix_form}.
\end{prop}

Based on \eqref{BPX_matrix_form} and \cref{prop:BPX_BCs}, we put forward the following parametrization
\begin{equation}
\label{modified_BPX}
    \widetilde{\mathcal{B}} = \sum_{k=1}^{L}\left(\widetilde{\alpha}_{k}\right)^2\widetilde{B}_{k}^{L}\widetilde{B}_{L}^{k},~\widetilde{B}_{l}^{L} = I_{l}\otimes\widetilde{\eta}_{L-l} + S_{l}\otimes \widetilde{\xi}_{L-l},\left(\widetilde{\xi}_{L-l}\right)_{2^{l}} = 0,~\widetilde{\eta}_{0} = 1,~\widetilde{\alpha}_{L} = 1,
\end{equation}
where $\widetilde{\alpha}_{k}, \widetilde{\eta}_{L-k}$ and $\widetilde{\xi}_{L-k}$ are free parameters that correspond to $\omega$ in \cref{algorithm:rho_optimization} and \cref{algorithm:kappa_optimization}. Chosen parametrization differs from \eqref{BPX_matrix_form} in two respects. First, we use $\widetilde{\xi}_{L-k}$ in place of $\widetilde{\xi}_{L-k}-\widetilde{\eta}_{L-k}$. Since both $\widetilde{\eta}_{L-k}$ and $\widetilde{\xi}_{L-k}$ are free parameters, both options lead to the same family of preconditioners. Second, we use $\left(\widetilde{\alpha}_{k}\right)^2$ in place of $\widetilde{\alpha}_{k}$. This choice among with conditions $\widetilde{\eta}_{0} = 1$ and $\widetilde{\alpha}_{L} = 1$ guarantee that $\widetilde{\mathcal{B}}$ is positive definite regardless of the choice of other parameters. Indeed, $\widetilde{\mathcal{B}}$ has a form $I + \sum_{k=1}^{L-1}\left(\widetilde{\alpha}_{k}\right)^2 \left(B_{L}^{k}\right)^TB_{L}^{k}$, that is, the sum of positive definite and positive semidefinite matrices.
Because of that, conditions of \cref{prop:optimization} apply and we can use parametric family \eqref{modified_BPX} to optimize condition number with \cref{algorithm:rho_optimization} and \cref{algorithm:kappa_optimization}. The last condition $\left(\widetilde{\xi}_{L-l}\right)_{2^{l}} = 0$ ensures that basis functions on level $l$ have the same support as the ordinary tent functions.

\section{Experiments}
\label{section:Numerical_examples}
Here we present the results of the optimization for a set of test problems. First, we give an overview of model equations and the discretization used and then comment on the performance of optimized BPX preconditioners.

\subsection{Model equations}

\begin{figure}
    \centering
    \subfloat[][Linear FEM]{
        \raisebox{12ex}
        {\begin{tabular}{c||c|c|c||c|c|c||c|c|c}
     & \multicolumn{3}{|c||}{BPX} & \multicolumn{3}{|c||}{\shortstack{$\phi_i$ are fixed\\optimized BPX}} & \multicolumn{3}{|c}{optimized BPX}\\
    \hline
    $L$ & $\rho$ & $\kappa$& $N$ &  $\rho$ & $\kappa$ & $N$ &  $\rho$ & $\kappa$ & $N$\\
    \hline
    $3$ & $0.611$ & $4.138$ & $5$ & $0.483$ & $2.866$ & $4$ & $0.332$ & $1.994$ & $3$\\
    \hline
    $4$ & $0.696$ & $5.58$ & $7$ & $0.554$ & $3.484$ & $4$ & $0.357$ & $2.109$ & $3$\\
    \hline
    $5$ & $0.744$ & $6.81$ & $8$ & $0.599$ & $3.983$ & $5$ & $0.367$ & $2.159$ & $3$\\
    \hline
    $6$ & $0.774$ & $7.845$ & $9$ & $0.629$ & $4.389$ & $5$ & $0.37$ & $2.174$ & $3$\\
    \hline
    $7$ & $0.794$ & $8.718$ & $10$ & $0.651$ & $4.724$ & $6$ & $0.373$ & $2.19$ & $3$\\
    \hline
    $8$ & $0.809$ & $9.456$ & $11$ & $0.667$ & $5.003$ & $6$ & $0.377$ & $2.21$ & $3$\\
\end{tabular}}
        \label{fig:Poisson_1D:FEM}
    }
    \subfloat[][Basis function]{
        \includegraphics[scale=0.35]{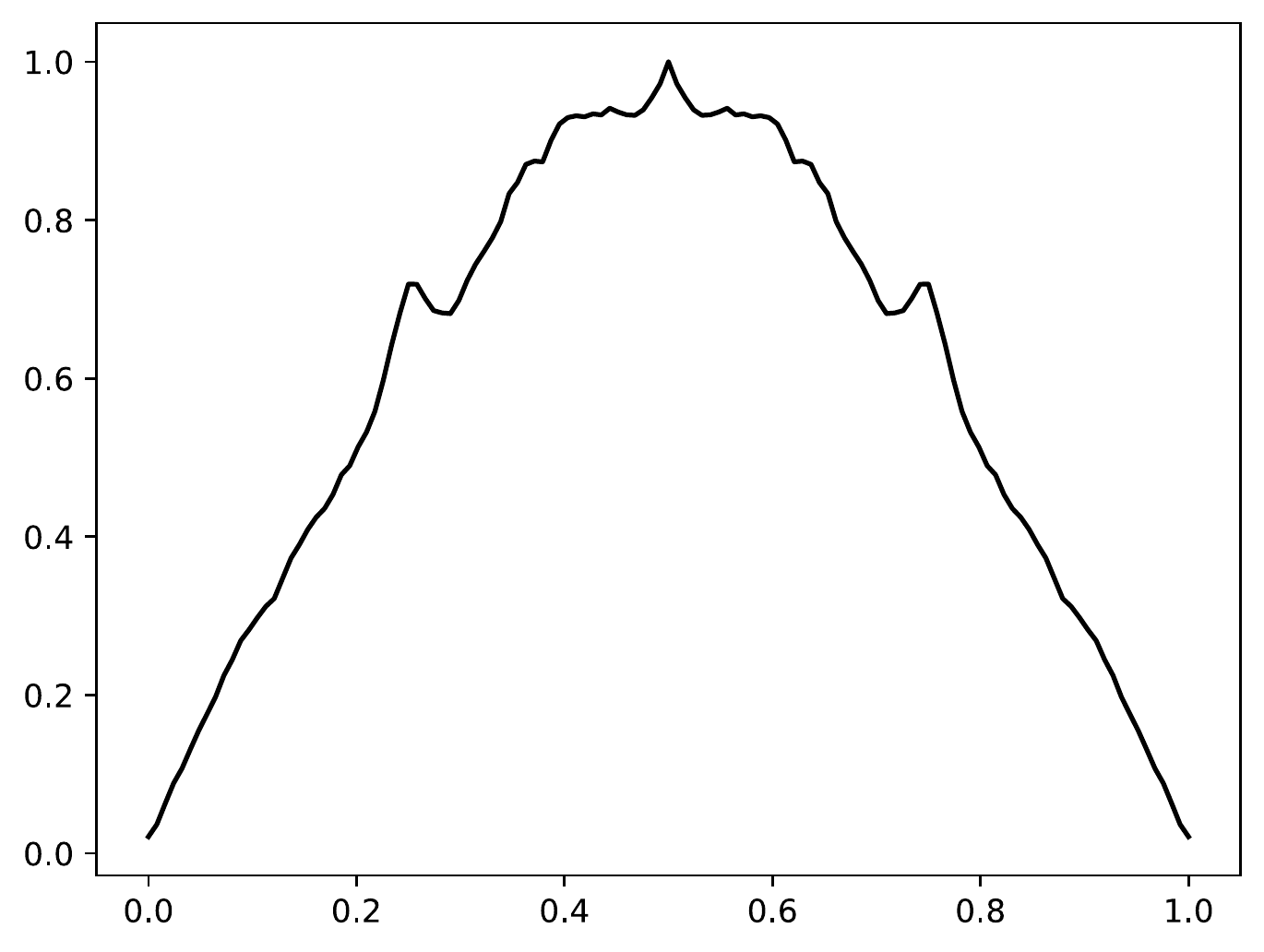}
        \label{fig:Poisson_1D:Basis}
    }
    \caption{Results of optimization and basis function for $1\text{D}$ Poisson equation \eqref{Poisson_equation}.}
    \label{fig:Poisson_1D}
\end{figure}

\subsubsection{Poisson equation}
Poisson equation appears in a variety of contexts, from continuum mechanics \cite[Sections 4.3, 5.1]{pletcher2012computational} to electrodynamics \cite[Section 1.7]{jackson1999classical}. It is also a standard test equation for multilevel solvers and preconditioners \cite[Section 1.4]{trottenberg2000multigrid}. The continuum boundary value problem reads
\begin{equation}
\label{Poisson_equation}
    -\frac{\partial^2 u(x, y)}{\partial x^2}-\frac{\partial^2 u(x, y)}{\partial y^2} = f(x),~x, y \in \left[0, 1\right]^2,~\left.u(x, y)\right|_{\partial\Gamma} = 0,
\end{equation}
here $\Gamma$ represents a domain, and $\partial\Gamma$ is a boundary. We use standard bilinear finite element discretization in $D=1$ and $D=2$ (see \cref{section:BPX_precs}), and also employ a high order compact scheme known as Mehrstellen \cite[Table VI]{collatz2012numerical}. Mehrstellen discretization corresponds to the stencil
\begin{equation}
    s = 
    \left[\begin{matrix}
        -1 & -4 & -1\\
        -4 & 20 & -4\\
        -1 & -4 & -1\\
    \end{matrix}\right],
\end{equation}
which can be used to construct a fourth and sixth-order accurate approximation to the Poisson equation if boundary conditions and right-hand side are sufficiently smooth \cite{rosser1975nine}.

\subsubsection{Helmholtz equation}
Helmholtz equation
\begin{equation}
\label{Helmholtz_equation}
    -\frac{\partial^2 u(x, y)}{\partial x^2}-\frac{\partial^2 u(x, y)}{\partial y^2} -k^2u(x,y)= f(x),~x, y \in \left[0, 1\right]^2,~\left.u(x, y)\right|_{\partial\Gamma} = 0,
\end{equation}
appears in the context of wave propagation problems \cite[Section 2.1]{erlangga2008advances}. For example, the Helmholtz equation needs to be solved at each time step in the semi-implicit discretization of governing equation of non-hydrostatic weather prediction models \cite[Section 4.1]{steppeler2003review}.

Because of the term $-k^2 u(x, y)$, bilinear finite element discretization can result in an indefinite matrix, especially for large $k$, which renders our method inapplicable. However, the value of $k$ can not be arbitrary on a given grid because of the pollution problem \cite{babuska1997pollution}. More precisely, unless $k^2h$ is sufficiently small, the solution to a discrete problem is of no use because it does not approximate an exact solution. Having this condition in mind, we choose $k$ small enough to have a positive definite problem.

\begin{figure}
    \centering
    \subfloat[][$k^2h=0.01$]{
        \begin{tabular}{c||c|c|c||c|c|c}
     & \multicolumn{3}{|c||}{BPX} & \multicolumn{3}{|c}{optimized BPX} \\
    \hline
    $L$ & $\rho$ & $\kappa$& $N$ &  $\rho$ & $\kappa$& $N$\\
    \hline
    $3$ & $0.621$ & $4.277$ & $5$ & $0.316$ & $1.922$ & $2$\\
    \hline
    $4$ & $0.701$ & $5.678$ & $7$ & $0.385$ & $2.254$ & $3$\\
    \hline
    $5$ & $0.746$ & $6.867$ & $8$ & $0.431$ & $2.515$ & $3$\\
    \hline
    $6$ & $0.774$ & $7.866$ & $10$ & $0.457$ & $2.685$ & $3$\\
\end{tabular}
        \label{fig:Helmholtz_Anisotropic:Helmholtz:a}
    }
    \subfloat[][$\epsilon=10$]{
        \begin{tabular}{c||c|c|c||c|c|c}
     & \multicolumn{3}{|c||}{BPX} & \multicolumn{3}{|c}{\shortstack{$s=1$\\optimized BPX}}\\
    \hline
    $L$ & $\rho$ & $\kappa$& $N$ &  $\rho$ & $\kappa$& $N$\\
    \hline
    $3$ & $0.919$ & $23.719$ & $28$ & $0.592$ & $3.9$ & $5$\\
    \hline
    $4$ & $0.956$ & $44.127$ & $51$ & $0.625$ & $4.339$ & $5$\\
    \hline
    $5$ & $0.967$ & $60.262$ & $70$ & $0.653$ & $4.766$ & $6$\\
    \hline
    $6$ & $0.973$ & $72.413$ & $84$ & $0.68$ & $5.25$ & $6$\\
\end{tabular}
        \label{fig:Helmholtz_Anisotropic:Anisotropic:a}
    }\\
    \subfloat[][$k^2h=0.1$]{
        \begin{tabular}{c||c|c|c||c|c|c}
     & \multicolumn{3}{|c||}{BPX} & \multicolumn{3}{|c}{optimized BPX}\\
    \hline
    $L$ & $\rho$ & $\kappa$& $N$ &  $\rho$ & $\kappa$& $N$\\
    \hline
    $3$ & $0.621$ & $4.277$ & $5$ & $0.316$ & $1.922$ & $2$\\
    \hline
    $4$ & $0.701$ & $5.678$ & $7$ & $0.385$ & $2.254$ & $3$ \\
    \hline
    $5$ & $0.746$ & $6.867$ & $8$ & $0.431$ & $2.515$ & $3$ \\
    \hline
    $6$ & $0.774$ & $7.866$ & $10$ & $0.457$ & $2.685$ & $3$ \\
\end{tabular}
        \label{fig:Helmholtz_Anisotropic:Helmholtz:b}
    }
    \subfloat[][$\epsilon=100$]{
        \begin{tabular}{c||c|c|c||c|c|c}
     & \multicolumn{3}{|c||}{BPX} & \multicolumn{3}{|c}{\shortstack{$s=2$\\optimized BPX}}\\
    \hline
    $L$ & $\rho$ & $\kappa$& $N$ &  $\rho$ & $\kappa$& $N$\\
    \hline
    $3$ & $0.974$ & $75.467$ & $87$ & $0.679$ & $5.235$ & $6$\\
    \hline
    $4$ & $0.991$ & $216.104$ & $249$ & $0.704$ & $5.763$ & $7$\\
    \hline
    $5$ & $0.996$ & $468.362$ & $540$ & $0.71$ & $5.9$ & $7$\\
    \hline
    $6$ & $0.997$ & $753.064$ & $867$ & $0.754$ & $7.145$ & $9$\\
\end{tabular}
        \label{fig:Helmholtz_Anisotropic:Anisotropic:b}
    }\\
    \subfloat[][$k^2h=1$]{
        \begin{tabular}{c||c|c|c||c|c|c}
     & \multicolumn{3}{|c||}{BPX} & \multicolumn{3}{|c}{optimized BPX}\\
    \hline
    $L$ & $\rho$ & $\kappa$& $N$ &  $\rho$ & $\kappa$& $N$\\
    \hline
    $3$ & $0.616$ & $4.213$ & $5$ & $0.317$ & $1.928$ & $3$\\
    \hline
    $4$ & $0.698$ & $5.612$ & $7$ & $0.387$ & $2.264$ & $3$ \\
    \hline
    $5$ & $0.744$ & $6.808$ & $8$ & $0.432$ & $2.519$ & $3$ \\
    \hline
    $6$ & $0.773$ & $7.817$ & $9$ & $0.457$ & $2.683$ & $3$ \\
\end{tabular}
        \label{fig:Helmholtz_Anisotropic:Helmholtz:c}
    }
    \subfloat[][$\epsilon=1000$]{
        \begin{tabular}{c||c|c|c||c|c|c}
     & \multicolumn{3}{|c||}{BPX} & \multicolumn{3}{|c}{\shortstack{$s=2$\\optimized BPX}}\\
    \hline
    $L$ & $\rho$ & $\kappa$& $N$ &  $\rho$ & $\kappa$& $N$\\
    \hline
    $3$ & $0.983$ & $118.948$ & $137$ & $0.694$ & $5.531$ & $7$\\
    \hline
    $4$ & $0.997$ & $578.133$ & $666$ & $0.735$ & $6.554$ & $8$\\
    \hline
    $5$ & $0.999$ & $1713.449$ & $1973$ & $0.763$ & $7.454$ & $9$\\
    \hline
    $6$ & $1.0$ & $4032.087$ & $4643$ & $0.807$ & $9.348$ & $11$\\
\end{tabular}
        \label{fig:Helmholtz_Anisotropic:Anisotropic:c}
    }
    \caption{First column (a, c, e) results for Helmholtz equation \eqref{Helmholtz_equation}, second column (b, d, f) results for anisotropic Poisson equation \eqref{Anisotropic_Poisson_equation}; $s$ refers to semicoarsening \eqref{modified_BPX_semicoarsening}.}
    \label{fig:Helmholtz_Anisotropic}
\end{figure}

\subsubsection{Anisotropic Poisson equation}
Anisotropic version of Poisson equation
\begin{equation}
\label{Anisotropic_Poisson_equation}
    -\frac{\partial^2 u(x, y)}{\partial x^2}-\epsilon\frac{\partial^2 u(x, y)}{\partial y^2} = f(x),~x, y \in \left[0, 1\right]^2,~\left.u(x, y)\right|_{\partial\Gamma} = 0,
\end{equation}
arises naturally in computational fluid dynamics when a refined or stretched grid is used to resolve a boundary layer, shock, or some other singularity \cite[Chapter 4]{liseikin2017grid}, \cite[Section 5.1.2]{trottenberg2000multigrid}. Parameter $\epsilon$ can also be related to the anisotropy of the physical system. For example, a crystal's permittivity can depend on the direction \cite[Chapter 9]{newnham2005properties}, so electrostatic boundary-value problems lead to an anisotropic Poisson equation.

\subsubsection{Biharmonic equation}
The only fourth-order equation we consider is biharmonic:
\begin{equation}
\label{Biharmonic_equation}
    \frac{\partial^4}{\partial x^4}u(x, y) + 2\frac{\partial^2}{\partial x\partial y}u(x, y) + \frac{\partial^4}{\partial y^4}u(x, y) = f(x, y),~x, y \in \left[0, 1\right]^2,~\left.u(x, y)\right|_{\partial\Gamma} = 0,~\left.\partial_{n} u(x, y)\right|_{\partial\Gamma} = 0,
\end{equation}
here $\partial_{n}$ is a derivative along the normal direction to the boundary $\partial \Gamma$. Applications of the Biharmonic equation include a description of fluid flows \cite{chen2004fast}, vibrating plates, Chladni figures \cite{gander2012chladni}, gravitation theory, and quantum mechanics  \cite[Introduction]{mak2018solving}. To discretize this equation, we use centered second-order finite difference approximation given by a $13$ point stencil
\begin{equation}
    s = 
    \left[\begin{matrix}
    &&1&&\\
    &2&-8&2&\\
    1&-8&20&-8&1\\
    &2&-8&2&\\
    &&1&&
    \end{matrix}\right],
\end{equation}
which should be modified appropriately near the boundaries \cite[Section 4]{tong1992multilevel} (see also \cite{gupta1979direct} and \cite{bramble1966second}).

\begin{figure}
    \centering
    \subfloat[][13-point stencil]{
        \raisebox{10ex}
        {\begin{tabular}{c||c|c|c||c|c|c}
     & \multicolumn{3}{|c||}{BPX} & \multicolumn{3}{|c}{optimized BPX}\\
    \hline
    $L$ & $\rho$ & $\kappa$& $N$ &  $\rho$ & $\kappa$& $N$\\
    \hline
    $3$ & $0.878$ & $15.367$ & $18$ & $0.846$ & $11.984$ & $14$\\
    \hline
    $4$ & $0.96$ & $48.717$ & $57$ & $0.878$ & $15.33$ & $18$ \\
    \hline
    $5$ & $0.988$ & $167.576$ & $193$ & $0.899$ & $18.9$ & $22$ \\
    \hline
    $6$ & $0.997$ & $617.095$ & $711$ & $0.945$ & $35.073$ & $41$\\
\end{tabular}}
        \label{fig:Biharmonic:FEM}
    }
    \subfloat[][Basis function]{
        \includegraphics[scale=0.35]{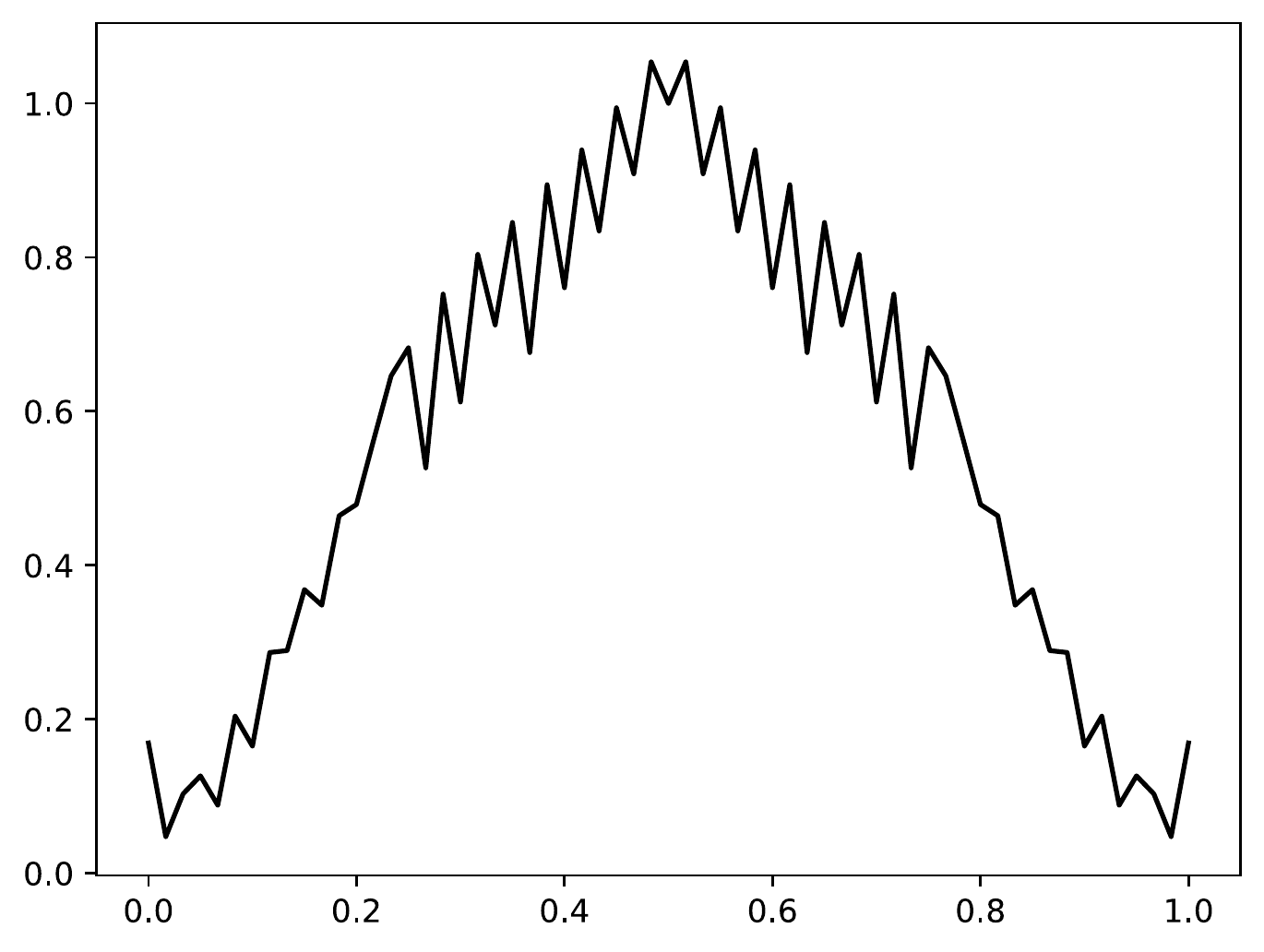}
        \label{fig:Biharmonic:Basis}
    }
    \caption{Results of optimization and basis function for the Biharmonic equation \eqref{Biharmonic_equation}.}
    \label{fig:Biharmonic}
\end{figure}

\subsubsection{Convection-diffusion equation}
When convective transport is present, the original diffusion equation needs to be modified as follows
\begin{equation}
\label{Convection_Diffusion_equation}
    -\frac{\partial^2 u(x, y)}{\partial x^2}-\frac{\partial^2 u(x, y)}{\partial y^2} + v_{x} u(x, y) + v_{y} u(x, y) = f(x, y),~x, y \in \left[0, 1\right]^2,~\left.u(x, y)\right|_{\partial\Gamma} = 0.
\end{equation}
The presence of $v_x$ and $v_y$ results in nonsymmetric matrix. This means \cref{prop:optimization} does not hold, but \cref{algorithm:rho_optimization} can be applied to optimize modified Richardson iteration. Since we employ bilinear finite element discretization (centered difference approximation), the stability restriction is given by Peclet condition $\max\left(|v_x|, |v_y|\right)\leq 2\big/h$.

\begin{figure}
    \centering
    \subfloat[][$v_x=-v_y=1\big/h$]{
        \begin{tabular}{c||c|c||c|c}
     & \multicolumn{2}{|c||}{BPX} & \multicolumn{2}{|c}{\shortstack{optimized \\BPX}}\\
    \hline
    $L$ & $\rho_3$ & $N$ & $\rho_3$ & $N$\\
    \hline
    $3$ & $0.629$ & $5$ & $0.398$ & $3$\\
    \hline
    $4$ & $0.741$ & $8$ & $0.554$ & $4$\\
    \hline
    $5$ & $0.797$ & $11$ & $0.649$ & $6$\\
    \hline
    $6$ & $0.829$ & $13$ & $0.690$ & $7$\\
\end{tabular}
    }
    \subfloat[][$v_x=-v_y=2\big/h$]{
        \begin{tabular}{c||c|c||c|c}
     & \multicolumn{2}{|c||}{BPX} & \multicolumn{2}{|c}{\shortstack{optimized \\BPX}}\\
    \hline
    $L$ & $\rho_3$ & $N$ & $\rho_3$ & $N$\\
    \hline
    $3$ & $0.787$ & $10$ &  $0.574$ & $5$\\
    \hline
    $4$ & $0.830$ & $13$ & $0.711$ & $7$\\
    \hline
    $5$ & $0.864$ & $16$ & $0.737$ & $8$\\
    \hline
    $6$ & $0.874$ & $18$ & $0.743$ & $8$\\
\end{tabular} 
    }
    \subfloat[][$v_x=-v_y=3\big/h$]{
        \begin{tabular}{c||c|c||c|c}
     & \multicolumn{2}{|c||}{BPX} & \multicolumn{2}{|c}{\shortstack{optimized \\BPX}}\\
    \hline
    $L$ & $\rho_3$ & $N$ & $\rho_3$ & $N$\\
    \hline
    $3$ & $0.855$ & $15$ &  $0.693$ & $7$\\
    \hline
    $4$ & $0.869$ & $17$ & $0.743$ & $8$\\
    \hline
    $5$ & $0.872$ & $17$ & $0.785$ & $10$\\
    \hline
    $6$ & $0.874$ & $18$ & $0.792$ & $10$\\
\end{tabular}
    }
    \caption{Results of optimization for convection-diffusion equation \eqref{Convection_Diffusion_equation}, $h$ is a distance between grid points on the finest grid. Note, that the value of a loss function \eqref{error_propagation_2} is listed, not an \textquote{exact} spectral radius.}
    \label{fig:Convection_Diffusion}
\end{figure}

\subsubsection{Diffusion with discontinuous coefficients}

In some situations, diffusion coefficient $a(x, y)$ in equation
\begin{equation}
\label{Disc_Diffusion_equation}
    -\frac{\partial}{\partial x}\left(a(x, y)\frac{\partial u(x, y)}{\partial x}\right)-\frac{\partial}{\partial y}\left(a(x, y)\frac{\partial u(x, y)}{\partial y}\right) = f(x, y),~x, y \in \left[0, 1\right]^2,~\left.u(x, y)\right|_{\partial\Gamma} = 0,
\end{equation}
is discontinuous along some curve or surface inside the computational domain. For example, this is the case in reservoir simulation \cite[Section 7.7.1]{trottenberg2000multigrid}, and the description of the neutron diffusion \cite{alcouffe1981multi}. For our experiments, we take 
\begin{equation}
\label{Disc_Diffusion_Coefficient}
    a(x, y) = g(x)+g(y),~g(x) = \sigma^{-1}{\sf Ind}\left[x <1\big/2\right] + \sigma{\sf Ind}\left[x \ge 1\big/2\right],
\end{equation}
where $\sigma$ is a parameter that controls the magnitude of the jump. The discretization we used is, again, FEM.

\subsubsection{Mixed derivative}
Another problem of interest is a Poisson equation with mixed derivative
\begin{equation}
\label{Mixed_equation}
    -\frac{\partial^2 u(x, y)}{\partial x^2}-\frac{\partial^2 u(x, y)}{\partial y^2} -2\tau\frac{\partial^2 u(x, y)}{\partial x\partial y}= f(x),~x, y \in \left[0, 1\right]^2,~\left.u(x, y)\right|_{\partial\Gamma} = 0.
\end{equation}
For $|\tau|>1$ the equation becomes hyperbolic, so it is interesting to look how optimization works for $\tau \simeq 1$.

\subsubsection{Implicit scheme for the heat equation}
The last equation that we consider comes from the trapezoidal discretization (in time) of the heat equation
\begin{equation}
\begin{split}
    &\frac{\partial u(x, y, t)}{\partial t} = \frac{\partial^2 u(x, y, t)}{\partial x^2} + \frac{\partial^2 u(x, y, t)}{\partial y^2},~x, y \in \left[0, 1\right]^2,~t\in\left[0,+\infty\right),\\ &\left.u(x, y, t)\right|_{t=0} = \phi(x, y),~\left.u(x, y)\right|_{\partial\Gamma} = 0.
\end{split}
\end{equation}
Let $A$ be a matrix that corresponds to a spatial FEM discretization of the right-hand side operator. It results in a system of ordinary differential equations
\begin{equation}
    \frac{d u_{i}(t)}{dt} = \sum_{j}A_{ij} u_{j}(t),~u_{i}(0) = \phi_{i}.
\end{equation}
Application of the trapezoidal rule leads to an unconditionally stable iteration
\begin{equation}
\label{Implicit_equation}
    \sum_{j}\left(I-\frac{\widetilde{\mu}}{2}A\right)_{ij}u_{j}^{n+1} = \sum_{j}\left(I+\frac{\widetilde{\mu}}{2}A\right)_{ij}u_{j}^{n}
\end{equation}
known as Crank-Nicolson scheme \cite[Section 16.4]{iserles2009first}. Here $\widetilde{\mu} = \Delta t$ is related to the Courant number $\mu = \Delta t/\Delta x^2$. Since matrix $\left(I-\frac{\widetilde{\mu}}{2}A\right)$ is symmetric positive definite for $\widetilde{\mu}\geq 0$ that needs to be inverted during each iteration, we test our preconditioner on this problem.

\subsection{Optimization results}
For all equations we use a symmetric form of both BPX \eqref{BPX_matrix_form} and modified BPX \eqref{modified_BPX} preconditioners. To access the results of optimization we list three related numbers: $\rho = \lambda_{\max}\left(I-\theta_{\text{opt}}BAB\right)$ -- a spectral radius of the optimal Richardson iteration for a given preconditioner, $\kappa = \lambda_{\max}(BAB)\big/\lambda_{\min}(BAB)$ -- spectral condition number, and $N$ -- the number of iteration needed to drop an error by $0.1$ with the optimal modified Richardson iteration in an arbitrary chosen norm, i.e., $\left\|e^{n+N}\right\|\big/\left\|e^{n}\right\| \leq 0.1$. The number of iterations $N$ is computed as $\left\lceil-1\big/\log_{10}\rho \right\rceil$, where $\left\lceil\cdot\right\rceil$ is the ceiling function. \footnote{This definition of $N$ guarantees $\left\|e^{n+N}\right\|\big/\left\|e^{n}\right\| \leq 0.1$ for normal iteration matrix $M(\omega, A)$. If $M(\omega, A)$ is not normal, $N$ holds as an estimation (see the discussion in Section~\ref{section:Direct_optimization_of_the_spectral_radius} after equation \eqref{linear_iteration}).}

In all cases, we use Dirichlet boundary conditions. Value of $L$ fixes the number of points along each direction to be $2^L-1$.

For all examples we employed \cref{algorithm:rho_optimization} with the loss function \eqref{error_propagation_2} ($N_\text{batch}=10$, $k=10$), ADAM optimizer \cite{kingma2014adam}, $N_{\text{epoch}}=500$, $N_{\text{inner}}=1$. Initial parameters $\widetilde{\alpha}, \widetilde{\eta}, \widetilde{\xi}$ of the modified BPX preconditioner \eqref{modified_BPX} were chosen such that the resulting matrix $\widetilde{\mathcal{B}}$ coincides with the BPX preconditioner \eqref{BPX_matrix_form}.

All algorithms were implemented in Julia \cite{bezanson2017julia} and available in a public repository \url{https://github.com/VLSF/neuralBPX}.

\begin{figure}
    \centering
    \subfloat[][$\sigma = 10$]{
        \begin{tabular}{c||c|c|c||c|c|c}
     & \multicolumn{3}{|c||}{BPX (r)} & \multicolumn{3}{|c}{optimized BPX (r)}\\
    \hline
    $L$ & $\rho$ & $\kappa$& $N$ &  $\rho$ & $\kappa$& $N$\\
    \hline
    $3$ & $0.727$ & $6.337$ & $8$ & $0.657$ & $4.834$ & $6$\\
    \hline
    $4$ & $0.898$ & $18.524$ & $22$ & $0.616$ & $4.213$ & $5$\\
    \hline
    $5$ & $0.964$ & $54.813$ & $64$ & $0.652$ & $4.746$ & $6$\\
    \hline
    $6$ & $0.986$ & $145.244$ & $168$ & $0.744$ & $6.811$ & $8$\\
\end{tabular}
    }
    \subfloat[][$\sigma = 100$]{
        \begin{tabular}{c||c|c|c||c|c|c}
     & \multicolumn{3}{|c||}{BPX (r)} & \multicolumn{3}{|c}{optimized BPX (r)}\\
    \hline
    $L$ & $\rho$ & $\kappa$& $N$ &  $\rho$ & $\kappa$& $N$\\
    \hline
    $3$ & $0.753$ & $7.082$ & $9$ & $0.614$ & $4.177$ & $5$\\
    \hline
    $4$ & $0.912$ & $21.831$ & $26$ & $0.692$ & $5.504$ & $7$\\
    \hline
    $5$ & $0.97$ & $66.626$ & $77$ & $0.739$ & $6.652$ & $8$\\
    \hline
    $6$ & $0.989$ & $186.977$ & $216$ & $0.809$ & $9.485$ & $11$\\
\end{tabular}
    }\\
    \subfloat[][$\tau = 0.5$]{
        \begin{tabular}{c||c|c|c||c|c|c}
     & \multicolumn{3}{|c||}{BPX} & \multicolumn{3}{|c}{optimized BPX}\\
    \hline
    $L$ & $\rho$ & $\kappa$& $N$ &  $\rho$ & $\kappa$& $N$\\
    \hline
    $3$ & $0.68$ & $5.255$ & $6$ & $0.45$ & $2.638$ & $3$\\
    \hline
    $4$ & $0.751$ & $7.044$ & $9$ & $0.511$ & $3.086$ & $4$ \\
    \hline
    $5$ & $0.79$ & $8.51$ & $10$ & $0.553$ & $3.479$ & $4$ \\
    \hline
    $6$ & $0.813$ & $9.685$ & $12$ & $0.577$ & $3.731$ & $5$ \\
\end{tabular}
    }
    \subfloat[][$\tau = 0.9$]{
        \begin{tabular}{c||c|c|c||c|c|c}
     & \multicolumn{3}{|c||}{BPX} & \multicolumn{3}{|c}{optimized BPX}\\
    \hline
    $L$ & $\rho$ & $\kappa$& $N$ &  $\rho$ & $\kappa$& $N$\\
    \hline
    $3$ & $0.817$ & $9.93$ & $12$ & $0.697$ & $5.599$ & $7$\\
    \hline
    $4$ & $0.89$ & $17.24$ & $20$ & $0.814$ & $9.781$ & $12$ \\
    \hline
    $5$ & $0.922$ & $24.787$ & $29$ & $0.864$ & $13.752$ & $16$ \\
    \hline
    $6$ & $0.935$ & $29.933$ & $35$ & $0.894$ & $17.811$ & $21$ \\
\end{tabular}
    }\\
    \subfloat[][$\widetilde{\mu}=h\big/2$]{
        \begin{tabular}{c||c|c|c||c|c|c}
     & \multicolumn{3}{|c||}{BPX} & \multicolumn{3}{|c}{optimized BPX}\\
    \hline
    $L$ & $\rho$ & $\kappa$& $N$ &  $\rho$ & $\kappa$& $N$\\
    \hline
    $3$ & $0.908$ & $20.723$ & $24$ & $0.067$ & $1.144$ & $1$\\
    \hline
    $4$ & $0.979$ & $94.196$ & $109$ & $0.033$ & $1.069$ & $1$ \\
    \hline
    $5$ & $0.995$ & $407.671$ & $470$ & $0.016$ & $1.033$ & $1$ \\
    \hline
    $6$ & $0.99$ & $1702.583$ & $1961$ & $0.017$ & $1.031$ & $1$ \\
\end{tabular}
    }
    \subfloat[][$\widetilde{\mu}=2\big/h$]{
        \begin{tabular}{c||c|c|c||c|c|c}
     & \multicolumn{3}{|c||}{BPX} & \multicolumn{3}{|c}{optimized BPX}\\
    \hline
    $L$ & $\rho$ & $\kappa$& $N$ &  $\rho$ & $\kappa$& $N$\\
    \hline
    $3$ & $0.641$ & $4.57$ & $6$ & $0.307$ & $1.885$ & $2$\\
    \hline
    $4$ & $0.729$ & $6.383$ & $8$ & $0.391$ & $2.287$ & $3$ \\
    \hline
    $5$ & $0.789$ & $8.47$ & $10$ & $0.505$ & $2.919$ & $4$ \\
    \hline
    $6$ & $0.845$ & $11.876$ & $14$ & $0.709$ & $5.875$ & $7$\\
\end{tabular}
    }
    \caption{Results of optimization for: first row (a, b) diffusion with discontinuous coefficients \eqref{Disc_Diffusion_Coefficient} ($(r)$ refers to rescaled version \eqref{modified_BPX_rescaled}), second row (c, d) Laplace operator with mixed derivative \eqref{Mixed_equation}, last row (e, f) matrix from Crank-Nicolson scheme \eqref{Implicit_equation}.}
    \label{fig:Disc_Mixed_Implicit}
\end{figure}

\subsubsection{Poisson equation}
We can see on \cref{fig:Poisson_2D} that for the $2\text{D}$ Poisson equation optimization successfully decreases the condition number. Moreover, it seems to grow slower compared to the original BPX preconditioner as the number of points increases \footnote{To estimate the growth rate we fit data using ordinary least squares with the model $\kappa(L) = c_1 + c_2 L$. For BPX preconditioner $\left(c_1, c_2\right) = \left(0.792, 1.196\right)$, and for the optimized BPX $\left(c_1, c_2\right) = \left(1.164, 0.264\right)$.}. To assess the contribution of the optimized basis functions, we perform additional optimization in $D=1$ with fixed basis functions. Results, given in \cref{fig:Poisson_1D}, indicate that optimization of the basis function leads to twice as small spectral radius compare to the situation when only scales are being optimized. The basis function itself is depicted in \cref{fig:Poisson_1D:Basis}. We can see that it is self-similar and seems to be well defined (in a sense that a subsampled basis function for $L_1>L_2$ is a good basis function for $L_2$). We can deduce that this function is a limit of some subdivision scheme \cite{rioul1992simple}, but we could not reliably define subdivision weights from our numerical experiments.
\subsubsection{Helmholtz equation}
The first column in \cref{fig:Helmholtz_Anisotropic} contains the results for Helmholtz equation \eqref{Helmholtz_equation} with $k^2h$ equal to $0.01$, $0.1$ and $1$. The results are similar to the one for the Poisson equation. However, if we further increase the number of points or $k$, the resulting matrix becomes indefinite, and the optimization breaks down. That means that with our approach, we cannot construct preconditioners for the Helmholtz equation. It is known that preconditioners for the Helmholtz equation significantly differ from preconditioners for Poisson-like equations (see \cite{erlangga2008advances} for the review), so this result is not surprising.
\subsubsection{Anisotrpoic Poisson equation}
The second column in \cref{fig:Helmholtz_Anisotropic} contains the results for anisotropic Poisson equation \eqref{Anisotropic_Poisson_equation} with $\epsilon$ equal to $10$, $100$ and $1000$. To cope with the anisotropy, we apply semicoarsening \cite[Section 5.1]{trottenberg2000multigrid}. Without semicoarsening a \textquote{projector} on the grid $M_{k}\times M_{k}$ ($M_{k}$ is as in \eqref{hierarchy_of_meshes}) reads $\widetilde{B}_{L}^{k}\otimes \widetilde{B}_{L}^{k}$. For semicoarsening the hierarchy of grids is modified, that is, in place of $M_{k}\times M_{k}$ we project on $M_{\min(k-s, 0)}\times M_{k}$, where $s$ quantifies the extent to which the grid along one direction is denser than a grid in the other direction. With this modification, a preconditioner itself takes a form
\begin{equation}
    \label{modified_BPX_semicoarsening}
    \widetilde{\mathcal{B}}_{s} = \sum_{k=1}^{L}\left(\widetilde{\alpha}_{k}\right)^2\left(\widetilde{B}_{\min(k-s, 0)}^{L}\otimes\widetilde{B}_{k}^{L}\right)\left(\widetilde{B}_{L}^{\min(k-s, 0)}\otimes\widetilde{B}_{L}^{k}\right).
\end{equation}
As a result, the coarsening is delayed for $y$ because $\epsilon>1$ in \eqref{Anisotropic_Poisson_equation}, i.e., $y$ is a direction of the strong coupling.
Note that in \cref{fig:Helmholtz_Anisotropic} we the compare \eqref{modified_BPX_semicoarsening} with original BPX preconditioner. If semicoarsening is applied to the BPX preconditioner, the weights $\alpha_{k}$ need to be modified. Original weights $\alpha_{k}$ combined with semicoarsening lead to worse performance. We can see that the optimization was able to fix the weights correctly. Moreover, comparing to semicoarsening applied in the context of filtering preconditioners \cite{tong1992multilevel} we were able to perform more aggressive coarsening, i.e., to decrease the number of floating-point operations.
\subsubsection{Biharmonic equation}
Results for the biharmonic equation are given in \cref{fig:Biharmonic}. We can see that the BPX preconditioner is relatively inefficient. It was able to substantially decrease the condition number compared to the original matrix (this condition number is not listed), but still, the condition number is large and grows like $\kappa_{L+1} \simeq 4 \kappa_{L}$. Condition number for the optimized BPX preconditioner is not only smaller but grows like $\kappa_{L+1} \simeq 2 \kappa_{L}$. The basis function on \cref{fig:Biharmonic:Basis} does not seem to be stable in this case. Authors in \cite{tong1992multilevel} were managed to obtain a better preconditioner for the biharmonic equation using larger filters. The same applies to the case of multigrid solvers, where orders of interpolation $n_{i}$ and restriction $n_{r}$ operators should fulfill $n_{i}+n_{r}> n_{l}$ \cite[Remark 2.7.1]{trottenberg2000multigrid}, where $n_{l}$ is the order of the linear operator ($4$ in the case of biharmonic equation). Given that, we can suggest that by increasing the basis function's support, one can achieve a better condition number. We will study this elsewhere.

\subsubsection{Convection-diffusion equation}
Convection-diffusion equation leads to a non-symmetric matrix. Because of this, we do not list spectral condition number in \cref{fig:Convection_Diffusion}. Here optimization results in about twice as efficient solver, but the improvement becomes less pronounced for larger convection coefficient values.
\subsubsection{Diffusion with discontinuous coefficients}
Because neither BPX nor modified BPX account for the variation of coefficients, we used a rescaled version of preconditioner
\begin{equation}
    \label{modified_BPX_rescaled}
    \widetilde{\mathcal{B}}_{r} = \sum_{k=1}^{L}\left(\widetilde{\alpha}_{k}\right)^2\widetilde{B}_{k}^{L} D\left(B_{L}^{k}AB_{k}^{L}\right)^{-1/2}\widetilde{B}_{L}^{k},
\end{equation}
where $D(\cdot)$ denotes the diagonal part of the matrix. For the original BPX preconditioner we again insert a diagonal part in-between \textquote{projectors} and use $\alpha_{k}$ as in \eqref{BPX_matrix_form}. Results are given in the first row of \cref{fig:Disc_Mixed_Implicit}. It is evident that it is enough to recover the correct scales $\widetilde{\alpha}_{k}$. This was achieved by optimization which produces a good preconditioner regardless of scale.

The other option would be to perform a Jacobi preconditioning step $A \rightarrow D(A)^{-1/2}A D(A)^{-1/2}$ as explained in \cite[discussion after equation (5.2)]{brezina1997robust} and (in relation to diffusion with discontinuous coefficients) in \cite[Section 3.1]{wathen2015preconditioning}. If this kind of rescaling is performed, BPX becomes a reasonable preconditioner, and optimization leads to results similar to the observed ones for the Poisson equation.
\subsubsection{Mixed derivative}
Results can be found in the second row of \cref{fig:Disc_Mixed_Implicit}. We can see that optimization is better for smaller values of $\tau$, but when $\tau$ becomes closer to one, optimization deteriorates. 

\subsubsection{Implicit scheme for heat equation}
Results are in the third row of \cref{fig:Disc_Mixed_Implicit}. We study problem \eqref{Implicit_equation} in two regimes. The first one corresponds to small time steps $\widetilde{\mu} = h\big/2$ used when the transient dynamic is of interest. In this case $I - \left(\widetilde{\mu}\big/2\right) A \simeq I$ so the preconditioner is not needed. As a result, BPX applied in a naive manner increases the condition number. The alternative solution would be to apply BPX preconditioner to the second matrix only, i.e., $I - \left(\widetilde{\mu}\big/2\right) BAB$, which solves this problem. However, the goal was to access the optimization, so we keep this experiment. In the other regime $\widetilde{\mu} = 2\big/h$ and one is interested in steady-state. In this situation, optimization again helps to decrease the spectral condition number. The last regime related to the elliptic equation with a linear source (different sign compare to the Helmholtz equation) for which a robust preconditioner was constructed in \cite{griebel1995tensor} with the help of a sophisticated subspace splitting technique.

\section{Conclusion}
In this article, we study the direct optimization of the spectral condition number. We derive two new loss functions, demonstrate how they are related to the spectral condition number, and show how stochastic optimization can be used to construct locally optimal preconditioners. We test our approach on a parametric family of modified BPX preconditioners. Optimization results show that for a large class of linear equations, automatic construction of reasonable preconditioners is possible. We want to emphasize that for many equations above, other more specialized preconditioners are available. There are also robust Schwarz preconditioners that are applicable for a broad class of second-order elliptic problems (see \cite{galvis2010domain}, \cite{efendiev2012robust}). The proposed approach differs from the previous attempts in three respects. First, described algorithms allow for a black-box construction of preconditioners, should a suitable parametrization is available. That means it is theoretically possible to apply the proposed approach in the algebraic setting as well. Second, the resulting preconditioner is locally optimal. The technique developed in \cite{efendiev2012robust} undoubtedly leads to a robust preconditioner. However, there is no guarantee that the resulting preconditioner is optimal. Since we are using stochastic gradient descent to directly optimize the spectral condition number of a preconditioner system, we can be sure that we achieve locally optimal preconditioner.\footnote{We can not guarantee global optimality within a giving family of preconditioners. The practical approach would be to use numerical continuation as explained in \cite{katrutsa2020black}.} Third, proposed algorithms can be potentially applied to a wider class of linear problems, f.e., different discretizations and higher-order equations. As a downside, our approach currently is not practically applicable for real problems because optimization includes thousands of matrix-vector products. However, it could be possible to transfer from optimization to learning, i.e., to construct a model that can be trained on small matrices and applied on larger matrices as it was done for the multigrid method \cite{greenfeld2019learning}. This is the focus of our current investigations.

\section{Acknowledgement}
The work was supported by Ministry of Science and Higher Education grant No. 075-10-2021-068.

\bibliography{refs.bib}
\end{document}